\documentclass[10pt,a4paper]{amsart}

\usepackage[latin1]{inputenc}
\usepackage{amsfonts}
\usepackage{amssymb}
\usepackage[leqno]{amsmath}
\usepackage{amsthm}
\usepackage{amscd}
\usepackage[all,import]{xy}
\usepackage[mathscr]{eucal}
\usepackage{times,verbatim}

\newtheoremstyle{theorem}{8pt\vfill}{8pt\vfill}{\itshape}{}{\bfseries}{.}{.5em}{}
\newtheoremstyle{paragraph}{8pt\vfill}{8pt\vfill}{}{}{\bfseries}{}{.5em}{}
\swapnumbers
\theoremstyle{theorem}                                       
\newtheorem{thm}{Theorem}[section]
\newtheorem*{questionA}{Problem A}
\newtheorem*{questionB}{Problem B}
\newtheorem*{questionC}{Problem C}

\newtheorem{prop}[thm]{Proposition}
\newtheorem{cor}[thm]{Corollary}
\newtheorem{lemma}[thm]{Lemma}

\theoremstyle{definition}
\newtheorem{df}[thm]{Definition}
\newtheorem{ex}[thm]{Example}
\newtheorem{rem}[thm]{Remark}

\theoremstyle{paragraph}                                     

\newcommand{\qbinom}[2]{\left[\begin{array}{c}\!\!#1\!\!\\ \!\!#2\!\!\end{array}\right]_q}

\DeclareMathOperator{\id}{id}
\DeclareMathOperator{\pr}{pr}
\DeclareMathOperator{\Hom}{Hom}

\DeclareMathOperator{\Spec}{Spec}
\DeclareMathOperator{\Gr}{Gr}      
\DeclareMathOperator{\SL}{SL}      
\DeclareMathOperator{\GL}{GL}

\DeclareMathOperator{\Sch}{Sch}      

\def\rk{{\rm{rk}\,}}       

\def\cC{\mathcal{C}}
\def\cD{\mathcal{D}}
\def\cF{\mathcal{F}}
\def\cG{\mathcal{G}}

\def\cP{\mathcal{P}}
\def\cX{\mathcal{X}}
\def\cY{\mathcal{Y}}

\def\cG{\mathcal{G}}
\def\cO{\mathcal{O}}
\def\cQ{\mathcal{Q}}
\def\cT{\mathcal{T}}
\def\cN{\mathcal{N}}
\def\Mo{\mathfrak{Mo}}
\def\fp{\mathfrak{p}}
\def\N{\mathbb{N}}
\def\Z{\mathbb{Z}}

\def\G{\mathbb{G}}
\def\cO{\mathcal{O}}
\def\F{\mathbb{F}}
\def\Fun{{\F_1}}                             
\def\Ga{\mathbb{G}_a}
\def\Gm{\mathbb{G}_m}
\def\A{\mathbb{A}}
\def\P{\mathbb{P}}
\def\str{\mathrm{str}}
\def\weak{\mathrm{weak}}
\def\nat{\mathrm{nat}}

\title{Algebraic groups over the field with one element}
\author{Oliver Lorscheid}
\address{Max-Planck-Institut f\"ur Mathematik\\
Vivatsga\ss{}e, 7. D-53111, Bonn, Germany}
\email{oliver@mpim-bonn.mpg.de}

\begin{document}

\maketitle

\begin{abstract}
Remarks in a paper by Jacques Tits from 1956 led to a philosophy how a theory of split reductive groups over $\F_1$, the so-called field with one element, should look like. Namely, every split reductive group over $\Z$ should descend to $\F_1$, and its group of $\F_1$-rational points should be its Weyl group. We connect the notion of a torified variety to the notion of $\F_1$-schemes as introduced by Connes and Consani. This yields models of toric varieties, Schubert varieties and split reductive groups as $\Fun$-schemes. We endow the class of $\F_1$-schemes with two classes of morphisms, one leading to a satisfying notion of $\F_1$-rational points, the other leading to the notion of an algebraic group over $\F_1$ such that every split reductive group is defined as an algebraic group over $\F_1$.  Furthermore, we show that certain combinatorics that are expected from parabolic subgroups of $\GL(n)$ and Grassmann varieties are realized in this theory.
\end{abstract}

\section*{Introduction}

\noindent
The development of $\F_1$-geometry plays a key r\^ole in the program of translating Weil's proof of the Riemann hypothesis as shaped by Kurokawa (\cite{Kurokawa92}), Deninger (\cite{Deninger91}, \cite{Deninger92}, \cite{Deninger94}), Manin (\cite{Manin95}) and others in the early 1990s. But the first mention of the ``field with one element'' appeared in Jacques Tits' paper \cite{Tits56} from 1956 and his ideas are a main inspiration in the development of $\F_1$-geometry. Tits' remarks gave rise to a philosophy of groups and group actions over $\F_1$, which was first seriously treated by Connes and Consani in \cite{CC08}. For a further discussion of their results, see \cite[section 6.1]{LL09}. We will give an idea of this philosophy in the present introduction and show how to realize it in the following sections. 

While there are now general different frameworks for $\F_1$-geometry, a common theme is that $\F_1$ should be an object lying below the integers. this means that an $\F_1$-geometry should be a category $\Sch_\Fun$  with a terminal object $\ast_\Fun=\Spec_\Fun\F_1$ and a base extension functor $-\otimes_\Fun\Z$ from $\Sch_\Fun$ to the category $\Sch_\Z$ of schemes such that $\ast_\Fun\otimes_\Fun\Z$ is isomorphic to $\ast_\Z=\Spec\Z$. Given a candidate for $\Sch_\Fun$, it is natural to ask: which schemes have a model over $\F_1$, i.e. for which schemes $X$ does exist an object $\cX$ in $\Sch_\Fun$ such that $\cX_\Z:=\cX\otimes_\Fun\Z$ is isomorphic to $X$?

The viewpoint originating from Tits' paper is the following. A wide class of schemes of finite type over $\Z$ admit a polynomial $N(q)$ with integer coefficients as a \emph{counting function}, that is, $N(q)$ equals the number of $\F_q$-points of the scheme for every prime power $q$. First examples include affine spaces, projective spaces and Grassmannians:
\begin{align*}
 \#\A^n(\F_q) &= q^n, 
 & \#\P^{n-1}(\F_q) &=[n]_q,  
 &\text{and} &
 & \#\Gr(k,n)(\F_q) &= \qbinom nk
\end{align*}
where $[n]_q=q^{n-1}+\dotsb+q+1$ is the Gauss number, $[n]_q!=\prod_{i=1}^n[i]_q$ is the Gauss factorial and $\left[\substack{n\\k}\right]_q =\frac{[n]_q!}{[k]_q![n-k]_q!}$ is the Gauss binomial. Evaluating these polynomials at $q=1$ leads to interesting numbers, which should be thought of as the number $\#X(\F_1)$ of ``$\F_1$-rational points'' of the scheme $X$. Comparing cardinalities, we see that
\begin{align*}
 \#\A^n(\F_1) &= 1=\#\ast, 
 & \#\P^{n-1}(\F_1) &= n=\# M_n,
 & \#\Gr(k,n)(\F_1) &= \binom nk =\#M_{k,n}
\end{align*}
where $\ast$ is the one point set, $M_n=\{0,\dotsc,n-1\}$ and $M_{k,n}$ is the set of subsets of cardinality $k$ in $M_n$. We formulate a first problem.

\begin{questionA} \label{questionA}
 We seek a category $\Sch_\Fun$ with a terminal object $\ast_\Fun$ and a functor $-\otimes_\Fun\Z:\Sch_\Fun\to\Sch_\Z$ that contains objects $\A^n_\Fun$, $\P^{n-1}_\Fun$ and $\Gr(k,n)_\Fun$ (for $n\geq1$ and $0\leq k\leq n$) such that
 \begin{align*}
  \A_\Fun^n\otimes_\Fun\Z&\simeq\A^n&& \text{and} &\A_\Fun^n(\F_1)&\simeq \ast,\\
  \P_\Fun^{n-1}\otimes_\Fun\Z&\simeq\P^{n-1}&& \text{and} & \P_\Fun^{n-1}(\F_1)&\simeq M_n,\\
  \Gr(k,n)_\Fun\otimes_\Fun\Z&\simeq\Gr(k,n)&&\text{and} &\Gr(k,n)_\Fun(\F_1)&\simeq M_{k,n}\;.
 \end{align*}
\end{questionA}

There are already several approaches that give partial solutions to this problem. All suggestions for $\Sch_\Fun$ in literature contain models of toric varieties, which include $\A^n$ and $\P^{n-1}$. As we will see in the course of this text, $\Gr(k,n)$ has a model in the notion of $\F_1$-scheme as suggested by Connes and Consani in \cite{CC09}. However, in most categories, the set $\Hom(\Spec_\Fun\F_1,\cX)$ is not equal to what we expect as $\cX(\F_1)$. Note that $\P^{n-1}=\Gr(1,n)$, so part (i) follows from part (iii) of Problem A. In the present paper, we will introduce morphisms between $\F_1$-schemes as defined in \cite{CC09} and show that Problem A can be solved.  

Another interesting source of examples are split reductive groups $G$ with maximal split torus $T\simeq\Gm^r$, where $r$ is the rank of $G$. Let $N$ be the normalizer of $T$ in $G$ and $W=N(\Z)/T(\Z)$ the Weyl group of $G$. Let $B$ be a Borel subgroup of $G$ containing $T$. The Bruhat decomposition of $G$ (w.r.t. $T$ and $B$) is the natural morphism
$$ \coprod_{w\in W} BwB \longrightarrow G. $$
This morphism induces a bijection $\coprod BwB(k)\simeq G(k)$ of $k$-rational points for every field $k$. Since $BwB\simeq \Gm^r\times\A^{d_w}$ for certain $d_w\geq0$, the Bruhat decomposition shows that $G$ admits a polynomial counting function
$$ \# G(\F_q) \ = \ \sum_{w\in W} (q-1)^r q^{d_w}. $$

However, if the rank $r$ of $G$ is positive, then the value of this polynomial at $q=1$ is zero. A more interesting number of the counting polynomial $N(q)$ is
$$ \lim_{q\to1} \frac{N(q)}{(q-1)^\rho} $$
where $\rho$ is the order of vanishing of $N(q)$ in $q=1$, i.e. the lowest non-vanishing coefficient in the development of $N(q)$ in $q-1$. Note that in the previous cases of $\A^n$, $\P^{n-1}$ and $\Gr(k,n)$, we have $\rho=0$ and nothing changes regarding Problem A. In the case of a split reductive group $G$, we have $\rho=r$ and 
\begin{equation}\label{grouplim} \lim_{q\to1}\frac{\sum_{w\in W}(q-1)^rq^{d_w}}{(q-1)^r} \ = \ \lim_{q\to1} \sum_{w\in W} q^{d_w} \ = \ \# W. \end{equation}
It was indeed Tits' suggestion to interpret the Weyl group of a split reductive group as its set of $\F_1$-points. In the framework as above, this means that we should ask for a concept of ``algebraic groups over $\F_1$'' such that split reductive groups are defined as algebraic groups over $\F_1$ and such that their $\F_1$-points are isomorphic to their Weyl group. More precisely, consider the following problem. 
\begin{questionB}
 We seek a category $\Sch_\Fun$ with finite products and a terminal object $\ast_\Fun$ together with a functor $-\otimes_\Fun\Z:\Sch_\Fun\to\Sch_\Z$ that respects finite products and the terminal object such that for every split reductive group $G$ with group law $m:G\times G\to G$, there is a group object $\cG$ in $\Sch_\Fun$ with group law $\mu:\cG\times\cG\to \cG$ (in $\Sch_\Fun$) satisfying the following properties.
 \begin{enumerate}
  \item \label{B1} $\cG_\Z\simeq G$ as algebraic groups, i.e. there is an isomorphism $\varphi:\cG_\Z\to G$ such that
        $$ \xymatrix{\cG_\Z\times\cG_\Z\ar[rr]^{\mu_\Z}\ar[d]^{(\varphi,\varphi)} && \cG_\Z\ar[d]^\varphi \\ G\times G\ar[rr]^m && G} $$
		  commutes.
  \item \label{B2} $\cG(\F_1):=\Hom_\Fun(\ast_\Fun,\cG)$ together with the induced group structure is isomorphic to $W$ as a group such that the limit in \eqref{grouplim} is respected, i.e. the morphism
 $$ \begin{array}{cccccccc}\hspace{0pt}\sigma: & N(\Z)/T(\Z) \ = \ W & \stackrel\sim\longrightarrow & \cG(\F_1) & \stackrel{-\otimes_\Fun\Z}\longrightarrow& \cG_\Z(\Z)&\stackrel\varphi\longrightarrow & G(\Z) \\
	                                             &&&  (\ast_\Fun\to\cG)& \longmapsto & (\ast_\Z\to\cG_\Z) \end{array} $$
        maps each coset $nT(\Z)$ in $N(\Z)/T(\Z)$ to an element of $nT(\Z)\subset G(\Z)$.
 \end{enumerate}
\end{questionB}

For the following reason, this problem cannot be solved in general. Note that $\sigma$ is a group homomorphism, since it is a composition of group homomorphisms (for the fact that the base extension is a group homomorphism, see Proposition \ref{cart_functor}). Then the fact that $\sigma(nT(\Z))\subset nT(\Z)$, shows that $\sigma$ splits the short exact sequence of groups
 \begin{equation}\label{seq} \xymatrix@C=4pc{1\ar[r]&T(\Z)\ar[r]&N(\Z)\ar[r]&W\ar[r]\ar@{-->}@/_1pc/[l]_\sigma&1.} \end{equation}
 This, however, is not possible for every split algebraic group as the example of $\SL(2)$ shows.

In \cite{CC08}, Connes and Consani circumvent the lifting problem by using Tits' construction from \cite{Tits66}, which shows that a certain extension of \eqref{seq} splits for every split reductive group. In this way, the normalizer $N$ becomes a group object that is defined over ``$\F_{1^2}$'', but the group law of $G$ fails to be defined over $\F_{1^2}$ in general, cf.\ \cite[section 6.1]{LL09}. In this text, we will use a different method that allows us to define every split reductive group as a group object defined over $\F_1$ such that it has the expected group of $\F_1$-points. Namely, we will use the framework of $\F_1$-scheme given by Connes and Consani in \cite{CC09}, and introduce two different classes of morphisms between $\F_1$-schemes, one leading to a satisfying notion of $\F_1$-rational points, the other allowing models of all split reductive groups over $\F_1$.

Once we have established split reductive groups as group objects over $\F_1$, we can investigate group actions and ask whether a quotient exists. In the case of a standard parabolic subgroup $P$ of $\GL(n)$ of type $(k,n-k)$ acting on $\GL(n)$ by multiplication from the left, the quotient is $\Gr(k,n)$. Since $P$ is isomorphic to $\GL(k)\times\GL(n-k)\times\A^{k(n-k)}$ as a variety, it has a polynomial counting function, namely $N_P(q)=q^{k(n-k)}\cdot N_{\GL(k)}(q)\cdot N_{\GL(n-k)}(q)$, where $N_{\GL(r)}(q)=\sum_{w\in W}(q-1)^rq^{d_w}$ is the counting function of $\GL(r)$, where the Weyl group $W$ of $\GL(r)$ is isomorphic to $S_r$ for $r\in\{k,n-k\}$. The order of vanishing of $N_P(q)$ at $q=1$ is $k+(n-k)=n$ and the number of $\F_1$-rational points is
\begin{eqnarray*} 
\# P(\F_1) &=& \lim_{q\to1} \frac{q^{k(n-k)}\bigl(\sum_{\alpha\in S_k} (q-1)^kq^{d_\alpha}\bigr)\bigl(\sum_{\beta\in S_{n-k}} (q-1)^{n-k}q^{d_\beta}\bigr)}{(q-1)^n} \\ &=& \lim_{q\to1} q^{k(n-k)}\bigl(\sum_{\alpha\in S_k} q^{d_\alpha}\bigr)\bigl(\sum_{\beta\in S_{n-k}} q^{d_\beta}\bigr)\hspace{1cm} =\hspace{1cm} \#(S_k\times S_{n-k}) 
\end{eqnarray*}
The quotient of the action $l: P\times G\to G$ is the Grassmannian $\Gr(k,n-k)$ and the quotient of the action $l': (S_k\times S_{n-k})\times S_n \to S_n$ is $M_{k,n}$. There is a natural action $t:\GL(n)\times \Gr(k,n)\to\Gr(k,n)$ and a natural action $t':S_n\times M_{k,n}\to M_{k,n}$. This leads to the following problem.
\begin{questionC}
We seek a category $\Sch_\Fun$ with finite products and a terminal object $\ast_\Fun$ together with a functor $-\otimes_\Fun\Z:\Sch_\Fun\to\Sch_\Z$ that respects finite products and the terminal object such that there exist group objects $\cG$ and $\cP$, a group action $\lambda:\cP\times\cG\to\cG$ and a quotient $\cQ$ of $\lambda$ with the following properties.
 \begin{enumerate}
    \item There are isomorphisms $f: \cP_\Z\simeq P$ and $g: \cG_\Z\simeq \GL(n)$ of algebraic groups such that $\lambda_\Z:\cP_\Z\times \cG_\Z\to \cG_\Z$ is compatible with $l:P\times \GL(n)\to \GL(n)$, i.e.\ 
	 $$\xymatrix{\cP_\Z\times\cG_\Z\ar[d]_{(f,g)}\ar[rr]^{\lambda_\Z} && \cG_\Z\ar[d]^g\\ P\times\GL(n)\ar[rr]^l && \GL(n)} $$
	 commutes. There are bijections $\cP(\F_1)\simeq S_k\times S_{n-k}$ and $\cG(\F_1)\simeq S_n$ such that $\lambda(\F_1):\cP(\Fun)\times\cG(\Fun)\to\cG(\Fun)$ is compatible with $l':(S_k\times S_{n-k})\times S_n \to S_n$.
	 \item Let $\tau:\cG\times\cQ\to\cQ$ be the natural action on the quotient. There is an isomorphism $\cQ_\Z\simeq\Gr(k,n)$ of varieties such that $\tau_\Z:\cG_\Z\times\cQ_\Z\to\cQ_\Z$ is compatible with $t:\GL(n)\times \Gr(k,n)\to\Gr(k,n)$. There is a bijection $\cQ(\Fun)\simeq M_{k,n}$ such that $\tau(\Fun):\cG(\Fun)\times\cQ(\Fun)\to\cQ(\Fun)$ is compatible with $t':S_n\times M_{k,n}\to M_{k,n}$.
 \end{enumerate}
\end{questionC}

We will show that Problem C can be solved within the framework of this paper.

The text is organized as follows. In section \ref{group_objects}, we recall the basic facts about group objects in an arbitrary category with finite products and a terminal object. In section \ref{f1-schemes}, we introduce the notion of an $\F_1$-scheme as defined by Connes and Consani in \cite{CC09} and show that toric varieties descend to $\F_1$. In section \ref{torified_varieties}, we introduce the notion of a torified variety as defined by L\'opez Pe\`na and the author in \cite{LL09}. The important property is that every torified variety descends to $\F_1$. We recall from \cite{LL09} that toric varieties, Schubert varieties and split reductive groups are torified varieties and are thus defined over $\F_1$. 

In section \ref{strong_morphisms}, we define the notion of a strong morphism between $\F_1$-schemes. With relation to this class of morphisms, the sets $\cX(F_1)=\Hom_\Fun^\str(\Spec_\Fun\Fun, \cX)$ return for $\cX$ being the $\Fun$-schemes from Problems A and B the expected sets of $\Fun$-points. In particular, we solve Problem A. In section \ref{weak_morphisms}, we define the notion of a weak morphism between $\Fun$-schemes. In section \ref{cartesian_categories}, we introduce certain functors that allow us to pass group objects from one category to another. 

In section \ref{algebraic_f1_groups}, we define the notion of a group scheme over $\Fun$ as a group object in the category of $\Fun$-schemes together with weak morphisms. An algebraic group over $\Fun$ is a group scheme over $\Fun$ whose base extension to $\Z$ is an algebraic group. We show that extensions of finite groups by split tori, split reductive groups and successive extensions of the additive group scheme $\Ga$ descend to algebraic groups over $\Fun$. In particular, this solves a slight modification of Problem B. In section \ref{parabolic_subgroups}, we show that parabolic subgroups of $\GL(n)$ can be defined as algebraic groups over $\Fun$. We solve Problem C.

\textbf{Acknowledgements.} The author thanks the Max Planck Institute for the inspiring working environment. He thanks the organizers and the participants of the Nashville conference on $\F_1$ from May 2009 for many interesting discussions. He thanks Javier L\`opez Pe\~na and Lisa Carbone for stimulating conservations. He thanks Ethan Cotterill for his help with preparing the paper.

\section{Preliminaries on group objects}
\label{group_objects}

\noindent
To begin with, we review the concept of a group object and provide some facts that we will use later on. For more details, cf.\ \cite[Expos\'e 1, section 2]{SGA3I}, \cite[Section III.6]{MacLane98} and \cite[\S 0.1]{Mumford65}.

In this text, we say that a category $\cC$ is \emph{cartesian} if it contains finite products and a terminal object $\ast_\cC$. A cartesian category $\cC$ comes with the following canonical morphisms for all objects $A$ and $B$:
 \begin{itemize}
  \item an isomorphism $\pr_1: A\times\ast_\cC \to A$,
  \item the \emph{diagonal} $\Delta: A\to A\times A$,
  \item an isomorphism $\chi:A\times B\to B\times A$.
 \end{itemize}

A \emph{group object} in a cartesian category $\cC$ is a pair $(G,m)$, where $G$ is an object in $\cC$ and $m:G\times G\to G$ is a morphism in $\cC$ such that the multiplication 
$$ \begin{array}{cccc} m_\ast:& \Hom(X,G)\times\Hom(X,G)&\longrightarrow&\Hom(X,G)\\
                               &       (f,g)           &\longmapsto    & m\circ(f,g) \end{array} $$        
defines a group structure of $\Hom(X,G)$ for every object $X$ in $\cC$. We refer to $G$ as the group object, when the context is clear, and refer to $m$ as its \emph{group law}.  

There is an alternative characterization of group objects. 

\begin{prop}\label{alt_group}
 Let $G$ be an object and $m:G\times G\to G$ be a morphism in a cartesian category $\cC$ with terminal object $\ast_\cC$. Then $(G,m)$ is a group object if and only if there are morphisms $\epsilon: \ast_\cC\to G$ and $\iota:G\to G$ such that the following diagrams commute:
 \begin{enumerate}
  \item (associativity)\label{ob1}\quad
        $\xymatrix{G\times G\times G\ar[rr]^{(\id,m)}\ar[d]_{(m,\id)}&& G\times G\ar[d]^m\\ G\times G\ar[rr]^m&& G }$,
  \item (left (right) unit)\label{ob2}\quad
        $\xymatrix{G\times \ast_\cC\ar[rr]^{(\id,\epsilon)}\ar[dr]_{\pr_1}&& G\times G\ar[dl]^m\\ & G }$ 
		  \quad \begin{minipage}[t]{3cm}(resp.\ with $m$ \\ replaced by $m\circ\chi$),\end{minipage}
  \item (left (right) inverse)\label{ob3}\quad
        $\xymatrix{G\ar[r]^{\Delta}\ar[d]& G\times G\ar[r]^{(\id,\iota)}& G\times G\ar[d]^m\\ \ast_\cC\ar[rr]^\epsilon&& G }$,
		  \quad \begin{minipage}[t]{3cm}(resp.\ with $m$ \\ replaced by $m\circ\chi$).\end{minipage}
 \end{enumerate}
 Moreover, if $(G,m)$ is a group object, then $\epsilon$ and $\iota$ are unique with the property that the diagrams $\eqref{ob1}$--$\eqref{ob3}$ commute. The unit element of $\Hom(X,G)$ for any object $X$ in $\cC$ is the morphism $X\to\ast_\cC\stackrel\epsilon\to G$.
\end{prop}

We refer to $\epsilon$ as the \emph{unit of $G$ }and to $\iota$ as the \emph{inversion of $G$}. We sometimes say that a quadruple $(G,m,\epsilon,\iota)$ is a group object, when we want to label the morphisms $\epsilon$ and $\iota$ related to a group object $(G,m)$ explicitly. 

\begin{proof}
This proposition is standard. We give only a brief outline.

Let $(G,m)$ be a group object. If $\epsilon:\ast_\cC\to G$ is the unit of the group $\Hom(\ast_\cC,G)$ and $\iota:G\to G$ is the inverse of $\id:G\to G$ in the group $\Hom(G,G)$, then the diagrams $\eqref{ob1}$--$\eqref{ob3}$ commute, and these choices for $\epsilon$ and $\iota$ are unique. 
 
Conversely, assume that there are morphisms $\epsilon$ and $\iota$ such that the diagrams $\eqref{ob1}$--$\eqref{ob3}$ commute. Then for every object $X$, the multiplication of $\Hom(X,G)$ is defined by $\varphi\cdot\psi:=m\circ(\varphi,\psi)$ for all $\varphi,\psi\in\Hom(X,G)$, the map $X\to\ast_\cC\stackrel\epsilon\to G$ is the unit and $\iota\circ\varphi$ is the inverse of $\varphi\in\Hom(X,G)$.
\end{proof}

A \emph{homomorphism of group objects} $(G_1,m_1)$ and $(G_2,m_2)$ is a morphism $\varphi:G_1\to G_2$ such that
$$\xymatrix{G_1\times G_1\ar[rr]^{m_1}\ar[d]_{(\varphi,\varphi)} && G_1\ar[d]^\varphi \\  G_2\times G_2\ar[rr]^{m_2} && G_2 }$$
commutes. If the context is clear, we will simply say that $\varphi: G_1\to G_2$ is a homomorphism of group objects.

We collect some standard facts.

\begin{lemma}\label{group_homo}
 Let $\varphi:G_1\to G_2$ be a homomorphism of group objects.
 \begin{enumerate}
  \item For every object $X$, the map $\varphi_\ast:\Hom(X,G_1)\to\Hom(X,G_2)$ is a group homomorphism.
  \item Let $\epsilon_i$ and $\iota_i$ be the unit resp.\ inversion of $G$ for $i=1,2$. Then the diagrams
   $$ \xymatrix@R=0,24pc{&& G_1\ar[dd]^\varphi\\\ast_\cC\ar[urr]^{\epsilon_1}\ar[drr]_{\epsilon_2}\\ && G_2} \hspace{1cm}\text{and}\hspace{1cm}
       \xymatrix@R=1,5pc{G_1\ar[rr]^{\iota_1}\ar[d]_{\varphi}&& G_1\ar[d]^\varphi\\ G_2\ar[rr]^{\iota_2}&& G_2 } $$
	commute.	 \qed
 \end{enumerate}
\end{lemma}

Let $Y$ be an object and $(G,m)$ be a group object in $\cC$. A \emph{group action of $G$ on $Y$} is a morphism $\theta:G\times Y\to Y$ such that for every object $X$ in $\cC$, the map $\theta_\ast:\Hom(X,G)\times\Hom(X,Y)\to\Hom(X,Y)$ is an action of the group $\Hom(X,G)$ on the set $\Hom(X,Y)$. There is an alternative definition.

\begin{prop}\label{alt_group_action}
 Let $Y$ be an object, $(G,m,\epsilon,\iota)$ a group object and $\theta: G\times Y\to Y$ a morphism in $\cC$. Then $\theta$ is a group action if and only if the diagrams 
  $$\xymatrix{G\times G\times Y\ar[rr]^{(\id,\theta)}\ar[d]_{(m,\id)}&& G\times Y\ar[d]^\theta\\ G\times Y\ar[rr]^\theta&& Y} \hspace{0,8cm}\text{and}\hspace{0,8cm} \xymatrix{\ast_\cC\times Y\ar[rr]^{(\epsilon,\id)}\ar[dr]_{\pr_2}&& G\times Y\ar[dl]^\theta\\ & Y } $$
	  commute.\qed
\end{prop}

Let $\theta:G\times Y\to Y$ be a group action. Then $Q$ together with a morphism $p:Y\to Q$ is a \emph{(categorical) quotient of $\theta$} if the diagram 
\begin{equation}\label{quotient_diagram} \xymatrix{G\times Y\ar[rr]^\theta\ar[d]_{\pr_2} && Y\ar[d]^p \\ Y\ar[rr]^p && Q} \end{equation}
commutes and if for all morphisms $f:Y\to Z$ such that diagram \eqref{quotient_diagram} with $p:Y\to Q$ replaced by $f:Y\to Z$ commutes, there is a unique morphism $\bar f:Q\to Z$ such that $f=\bar f\circ p$.

A \emph{subgroup} of a group object $(G,m,\epsilon,\iota)$ is a group object $(H,m',\epsilon',\iota')$ together with a homomorphism $H\to G$ of group objects that is a monomorphism in $\cC$. By Lemma \ref{group_homo}, we can think of $(m',\epsilon',\iota')$ as the restriction of $(m,\epsilon,\iota)$ to $H$ and when will suppress the formal difference between $(m',\epsilon',\iota')$ and $(m,\epsilon,\iota)$ in the notation, when the context is clear.

A subgroup $H$ of group object $(G,m,\epsilon,\iota)$ \emph{acts on a subset $Y$ of $G$ by conjugation} if the image of the morphism
$$ \xymatrix@C=3pc{\bar\theta: H\times Y\ar[r]^{(\Delta,\id)}& H\times H\times Y\ar[r]^{(\id,\chi)}& H\times Y\times H\ar[r]^{(m,\iota)}&G\times W\ar[r]^m& G}$$
is contained in $Y$, i.e.\ $\bar\theta$ factors through a group action $\theta:H\times Y\to Y$, which we call the \emph{conjugation of $H$ on $Y$}. A \emph{normal subgroup of $G$} is a subgroup $N$ of $G$ on which $G$ acts by conjugation.

\begin{lemma}\label{quotient_group}
 If $N$ is a normal subgroup of a group object $G$ and $Q$ is a quotient of the conjugation $G\times N\to N$, then $Q$ inherits a natural structure of a group and we call $Q$ the \emph{quotient group} of $G$ by $N$.
\end{lemma}

The \emph{direct product of group objects $(G_1,m_1)$ and $(G_2,m_2)$} is the product $G_1\times G_2$ together with the pair $m=(m_1,m_2)$ as group law, which is easily seen to define a group object. 

Let $(N,m_N)$ and $(H,m_H)$ be group objects and $\theta:H\times N\to N$ a group action that respects the group law $m_N$ of $N$, i.e.\ if we define the \emph{change of factors along $\theta$} as
$$ \xymatrix@C=3,5pc{\chi_\theta:\quad H\times N\ar[r]^{(\Delta,\id)}& H\times H\times N\ar[r]^{(\id,\theta)}&H\times N\ar[r]^\chi& N\times H,} $$
then the diagram 
$$ \xymatrix@C=4pc@R=0,8pc{H\times N\times N\ar[r]^{(\id,m_N)}\ar[dd]_{(\chi_\theta,\id)}&H\times N\ar[dr]^\theta\\ &&N\\ N\times H\times N\ar[r]^{(\id,\theta)}& N\times N\ar[ur]_{m_N}} $$
commutes. Then the morphism
$$\xymatrix@C=4pc{m_\theta:\quad N\times H\times N\times H\ar[r]^{\quad(\id,\chi_\theta,\id)}&N\times N\times H\times H\ar[r]^{\quad\quad(m_N,m_H)}&N\times H}$$
is a group law for $G=N\times H$. We say that $G$ is the \emph{semidirect product of $N$ with $H$ w.r.t.\ $\theta$} and write $G=N\rtimes_\theta H$. The group object $N$ is a normal subgroup of $G$ with quotient group $H$, and $H$ is a subgroup of $G$ that acts on $N$ by conjugation. The conjugation $H\times N\to N$ equals $\theta$. If $\theta: H\times N\to N$ is the canonical projection to the second factor of $H\times N$, then $N\rtimes_\theta H$ is equal to the direct product of $N$ and $H$ (as group object).

\begin{lemma}\label{semidirect}
 Let $G$ be the semidirect product of $N$ with $H$ w.r.t.\ $\theta$ and $G'$ be the semidirect product of $N'$ with $H'$ w.r.t.\ $\theta'$. If there are group homomorphisms $\varphi:N\to N'$ and $\psi:H\to H'$ such that
$$\xymatrix@C=4pc{N\times H\ar[rr]^\theta\ar[d]_{(\varphi,\psi)}&&N\ar[d]^\varphi\\ N'\times H'\ar[rr]^{\theta'}&&N'}$$
commutes, then $G=N\rtimes_\theta H\stackrel{(\varphi,\psi)}\longrightarrow N'\rtimes_{\theta'}H'=G'$ is a homomorphism of group objects.\qed 
\end{lemma}

Let $\cC$ and $\cD$ be cartesian categories. We say that a functor $\mathcal F:\cC\to\cD$ is \emph{cartesian} if $\cF(A\times B)\simeq\cF(A)\times\cF(B)$ and $\cF(\ast_{\cC})\simeq\ast_{\cD}$.

\begin{prop}\label{cart_functor}
 Let $\cF:\cC\to\cD$ be a cartesian functor between cartesian categories and $(G,m)$ be a group object in $\cC$. Then $(\cF(G), \cF(m))$ is a group object in $\cD$, and for every object $X$ in $\cC$, the map $\Hom_\cC(X,G)\to\Hom_\cD(\cF(X),\cF(G))$ sending $\varphi$ to $\cF(\varphi)$ is a group homomorphism.
\end{prop}

\begin{proof}
 By functoriality, $\cF$ maps the commutative diagrams $\eqref{ob1}$--$\eqref{ob3}$ to commutative diagrams. Since $\cF$ respects products and the terminal object, these diagrams verify that $(\cF(G), \cF(m))$ is a group object in $\cD$.
 
The last statement of the proposition follows from the equality
$$ \cF(\varphi\cdot\psi) = \cF(m\circ(\varphi,\psi)) = \cF(m)\circ\cF(\varphi,\psi) = \cF(m)\circ(\cF(\varphi),\cF(\psi)) = \cF(\varphi)\cdot\cF(\psi) $$
for any two morphisms $\varphi$ and $\psi$ in $\Hom(X,G)$.
\end{proof}

By functoriality, Lemma \ref{semidirect} and Proposition \ref{cart_functor} imply immediately the following.

\begin{cor}\label{cart_semidirect}
 Let $G=N\rtimes_\theta H$ be the semidirect product of $N$ with $H$ w.r.t.\ $\theta$ in $\cC$ and let $\cF:\cC\to\cD$ be a cartesian functor. Then $\cF(\cG)$ is the semidirect product $\cF(N)\rtimes_{\cF(\theta)}\cF(H)$ of $\cF(N)$ with $\cF(H)$ w.r.t.\ $\cF(\theta)$ in $\cD$.\qed
\end{cor}

In the following, a \emph{variety} means a reduced scheme of finite type over $\Z$. A \emph{group scheme} is a group object in the category $\Sch_\Z$ of schemes. An \emph{algebraic group} is group scheme that is a variety.

\section{Schemes over $\F_1$}
\label{f1-schemes}

\noindent
In this section, we review the definition of a scheme over $\F_1$ as given by Connes and Consani in \cite{CC09}. This notion combines the earlier ideas of \cite{CC08} and \cite{Soule04} with \cite{Deitmar05} and \cite{TV08}.

We begin with recalling the notion of an $\Mo$-scheme. For details, see \cite[section 3.6]{CC09}; also cf.\ \cite{Deitmar05}. Let $\Mo$ be the category of commutative monoids $M$ with $1$ (\emph{monoids} for short) and with $0$, i.e.\ an element satisfying $0\cdot a=0$ for all $a\in M$. A morphism between monoids with $0$ is a multiplicative map that sends $1$ to $1$ and $0$ to $0$. A \emph{monoidal space (with 0)} is a topological space $X$ together with a structure sheaf $\cO_X$ with values in $\Mo$. A \emph{morphism of monoidal spaces} is a continuous map together with a morphism of sheaves. Since direct limits exist in $\Mo$, it is possible to define stalks $\cO_{X,x}$ for every point $x\in X$. 

A \emph{prime ideal} of a monoid $M$ with $0$ is a subset $\fp$ of $M$ containing $0$ such that $\fp M\subset M$ and $M-\fp$ is a multiplicative subset containing $1$. It is possible to define localisations of $M$ at multiplicative subsets and to endow the set of prime ideals of $M$ with a Zariski topology in the same way as it is done for rings, since these constructions use only multiplicative structure. This defines a monoidal space $\Spec_\Mo M$, called the \emph{spectrum of $M$}. An \emph{$\Mo$-scheme} is a monoidal space that is locally isomorphic to the spectrum of a monoid. A morphism of $\Mo$-schemes is a morphism of monoidal spaces. Let $\Sch_\Mo$ denote the category of $\Mo$-schemes.

The category $\Sch_\Mo$ is cartesian. The terminal object is $\ast_\Mo=\Spec\{0,1\}$ and the product is locally given by $\Spec A\times\Spec B=\Spec (A\wedge B)$, where $A\wedge B$ denotes the smash product of $A$ and $B$ with respect to $0$ as base point.

There is a \emph{base extension functor} $\tilde X\mapsto \tilde X_\Z=\tilde X\otimes_\Fun\Z$
from $\Sch_\Mo$ to the category $\Sch_\Z$ of schemes (over $\Z$), which is locally described by
$$ \Spec_\Mo M \quad \longmapsto \quad \Spec\ \bigl(\Z[M]/(1\cdot0_M-0_{\Z[M]})\bigr), $$
where $\Z[M]$ is the semi-group ring of $M$, $0_M$ is the zero of $M$ and $0_{\Z[M]}$ is the zero of $\Z[M]$. This functor is cartesian.

A \emph{scheme over $\F_1$} (or \emph{$\F_1$-scheme}) is a triple $\cX=(\tilde X, X, e_X)$, where $\tilde X$ is an $\Mo$-scheme, $X$ is a scheme and $e_X: \tilde X_\Z\to X$ is a morphism of schemes (called the \emph{evaluation map}) such that $e_X(k):\tilde X_\Z(k)\to X(k)$ is a bijection for every field $k$. 
\begin{rem}
 An $\F_1$-scheme $\cX=(\tilde X, X, e_X)$ is \emph{locally of finite type} if $X$ is locally of finite type. For the sake of simplicity, we will assume for the rest of this text that \emph{all schemes over $\F_1$ are locally of finite type}.
\end{rem}


There is a natural choice of morphism of $\F_1$-schemes as a morphism between the underlying $\Mo$-schemes together with a morphism between the underlying schemes that are compatible with the evaluation maps. However, this notion of morphism is not suitable for a theory of algebraic groups over $\F_1$ as the only group laws that are of this nature are extension of finite groups by split tori (cf.\ Remark \ref{semidirect_remark}). We postpone the task to define the appropriate notion of morphism to a later section.

The \emph{base extension functor} $-\otimes_\Fun\Z$ associates to an $\F_1$-scheme $\cX=(\tilde X,X,e_X)$ the scheme $\cX_\Z:=X$.

\begin{ex}\label{ex:fun_from_mo}
To every $\Mo$-scheme $\tilde X$, we can associate the $\F_1$-scheme $\cX=(\tilde X,\tilde X_\Z, \id_{\tilde X_\Z})$. We have that $\tilde X\otimes_\Mo\Z=\cX\otimes_\Fun\Z$. We give first examples of $\F_1$-schemes of this kind. The affine line $\A^1_\Mo$ is the spectrum of the monoid $\{T^i\}_{i\in\N}\amalg\{0\}$ and the associated $\F_1$-scheme $(\A^1_\Mo,\A^1,\id_{\A^1})$ is a model of the affine line over $\F_1$. The multiplicative group $\G_{m,\Mo}$ is the spectrum of the monoid $\{T^i\}_{i\in\Z}\amalg\{0\}$, which defines $\G_{m,\Fun}=(\G_{m,\Mo},\Gm,\id_{\Gm})$ and base extends to $\Gm$ as desired. Both examples can be extended to define a model $(\A^n_\Mo,\A^n,\id_{\A^n})$ of the $n$-dimensional affine space over $\F_1$ and $\G_{m,\Fun}^n=(\G^n_{m,\Mo},\Gm^n,\id_{\Gm^n})$ by considering multiple variables $T_1,\ldots,T_n$. 

More generally, every toric variety $X$ gives rise to a connected $\Mo$-scheme $\tilde X$ such that $\tilde X_\Z\simeq X$. Thus toric varieties can be realized as $\F_1$-schemes. On the other hand, the only varieties that have models as $\Mo$-schemes are toric varieties. These observations are essentially due to Deitmar (\cite{Deitmar08}, consider also \cite[Thm.\ 4.1]{LL09}). 
\end{ex}

\section{Torified varieties}
\label{torified_varieties}

\noindent
We review the definition of a torified variety as introduced by Javier L\'opez Pe\~na and the author in \cite{LL09}. The connection to schemes over $\F_1$ is immediate and delivers a rich class of examples including Grassmannians and split reductive groups.

A \emph{torified scheme} is a scheme $X$ together with a \emph{torification} $e_X: T\to X$, that is, a morphism $T\to X$, where $T$ is a disjoint union $T=\coprod_{i\in I} \Gm^{d_i}$ of split tori, such that $e_X(k): T(k)\to X(k)$ is a bijection for all fields $k$. A \emph{torified variety} is a torified scheme $X$ that is a variety.

\begin{rem}\label{different_defs}
 The definition of a torified scheme given here differs from the original one given in \cite{LL09}. Namely, in \cite{LL09}, one meets the additional condition that the restrictions $e_X\vert_{\Gm^{d_i}}:\Gm^{d_i}\to X$ are immersions for all $i\in I$. The aim of \cite{LL09} was to establish examples of $\F_1$-varieties in the sense of the papers \cite{Soule04} and \cite{CC08}. For the aim of the present text, we do not need this additional property and thus work with the simplified (and more general) definition. 
 
The two definitions are close to each other, since every morphism $e: \Gm^d\to X$ from a split torus to a scheme $X$ is locally closed and injective. It is, however, not clear to me whether the morphisms $\cO_{X,e(y)}\to\cO_{\Gm^d,y}$ between the stalks are surjective for all $y\in\Gm^d$, which is the missing property for $e$ to be an immersion.
\end{rem}

Note that $T=\tilde X_T\otimes_\Mo\Z$ for $\tilde X_T=\coprod_{i\in I} \G_{m,\Mo}^{d_i}$. This yields immediately:

\begin{lemma}\label{fun_from_torified}
 Every torified scheme $X$ with torification $e_X:T\to X$ defines an $\F_1$-scheme $(\tilde X_T,X,e_X)$
\end{lemma}

In \cite[section 1.3]{LL09}, we find examples of torified varieties. We will recall these briefly. 

\begin{ex}[Toric varieties] \label{fun_from_toric}
The decomposition of a toric variety $X$ with torus action $T\times X\to X$ into the orbits of this action provides a torification of $X$. This establishes models of toric varieties as $\F_1$-schemes, again. 

We treat the example of a split torus and affine space in more detail. The split torus $\Gm^r$ has the trivial torification $\G^r_{m,\Fun}\to\G^r_{m,\Fun}$. With that, we obtain the same $\F_1$-scheme $\G^r_{m,\Fun}=(\G^r_{m,\Mo},\Gm^r,\id_{\Gm^r})$ as in Example \ref{ex:fun_from_mo}. The affine space $\A^d$ has a decomposition into tori $\Gm^I=\Spec\Z[T_i,T_i^{-1}]_{i\in I}$, where $I$ ranges through all subsets of $\{1,\dotsc,d\}$. The embedding $\Gm^I\hookrightarrow\A^d$ is given by the algebra homomorphism
$$ \begin{array} {ccc} \Z[T_1,\dots,T_d] & \longrightarrow & \Z[T_i,T_i^{-1}]_{i\in I}\;. \\
                                T_l      & \longmapsto     & \begin{cases}T_l&\text{if }l\in I\\0&\text{if }l\notin I\end{cases}  \end{array} $$
In particular, there is a unique torus of dimension $0$, which is embedded into the origin of $\A^d$. This defines the $\F_1$-scheme $\A^d_\Fun=(\tilde \A^d,\A^d,e_{\A^d})$. Note that the topological space of $\tilde\A^d$ is discrete, while the $\Mo$-scheme $\A^d_\Mo$ given by the canonical structure of $\A^d$ as toric variety (as described in Example \ref{ex:fun_from_mo}) is connected. Thus the $\F_1$-schemes $(\A^d_\Mo,\A^d,\id_{\A^d})$ and $(\tilde \A^d,\A^d,e_{\A^d})$ are not the same, but we will see in Remark \ref{models_of_toric} that they become isomorphic when we endow $\F_1$-schemes with ``strong morphisms''.
\end{ex}

\begin{ex}[Schubert varieties] \label{fun_from_schubert}
Another class of examples is Schubert varieties, which in particular includes Grassmann and flag varieties. Schubert varieties allow a decomposition into affine spaces that can be further decomposed into tori. In the case of the Grassmannian $\Gr(k,n)$, we have a Schubert decomposition
$$ \coprod_{w\in M_{k,n}} \A^{d_w} \longrightarrow \Gr(k,n), $$
which induces a bijection on $k$-points for every field $k$. The affine spaces $\A^{d_w}$ can be further decomposed into split tori, what yields a torification of $\Gr(k,n)$ and consequently a model $\Gr(k,n)_\Fun=(\tilde{\Gr}(k,n),\Gr(k,n),e_{\Gr(k,n)})$ of the Grassmannian over $\F_1$. Note that the $0$-dimensional tori in this torification stay in bijection with $M_{k,n}$.
\end{ex}

\begin{ex}[Split reductive groups]
\label{fun_from_split_red}
The last class of examples discussed in \cite{LL09} are split reductive groups $G$. For a definition, see \cite[Expos\'e XIX, Def.\ 2.7]{SGA3III}. Let $T\simeq\Gm^r$ be a maximal split torus of $G$, where $r$ is the rank of $G$. Let $N$ be the normalizer of $T$ in $G$ and $W=N(\Z)/T(\Z)$ the Weyl group of $G$. Let $B$ be a Borel subgroup of $G$ containing $T$. The Bruhat decomposition $\coprod_{w\in W} BwB\to G$ can be refined to a decomposition into split tori in the following way. (Note that we identify the coset $w\in W$ with the corresponding subvariety of $G$, which is isomorphic to $\Gm^r$). For every $w\in W$, we can choose an isomorphism $BwB\simeq\Gm^r\times\A^{d_w}$ for a certain $d_w\geq0$ as varieties. Therefore $BwB$ is toric and thus torified. More precisely, we can choose a torification of $BwB\simeq\Gm^r\times\A^{d_w}$ that contains $\Gm^r\simeq w\hookrightarrow BwB$ as $r$-dimensional torus; all other tori in the torification are of dimension larger than $r$. This provides a torification of $G$ that restricts to a torification of $N$ into $r$-dimensional tori, indexed by $W$.
\end{ex}

We collect the results obtained by these examples using Lemma \ref{fun_from_torified}.

\begin{prop}
 There are $\F_1$-schemes $\G^r_{m,\Fun}$, $\A^d_\Fun$ (Example \ref{fun_from_toric}), $\Gr(k,n)_\Fun$ (Example \ref{fun_from_schubert}) such that
 \begin{align*}
  \G^r_{m,\Fun}\otimes_\Fun\Z&\simeq\Gm^r,&
  \A^d_\Fun\otimes_\Fun\Z&\simeq\A^d,&
  \Gr(k,n)_\Fun\otimes_\Fun\Z&\simeq\Gr(k,n),
 \end{align*}
and there is a $\F_1$-scheme $\cG$ for every split reductive group $G$ (Example \ref{fun_from_split_red}) such that $\cG_\Z\simeq G$.
\end{prop}

\section{Strong morphisms}
\label{strong_morphisms}

\noindent
In this section, we define a class of morphisms between $\F_1$-schemes that produces the expected sets of $\F_1$-points for affine and projective space, Grassmann varieties and split reductive groups as formulated in Problems A and B of the introduction.

\noindent
Let $\cX=(\tilde X,X,e_X)$ and $\cY=(\tilde Y,Y,e_Y)$ be $\F_1$-schemes. Then we define the \emph{rank of a point $x$} of the underlying topological space $\tilde X$ as $\rk x:=\rk\cO_{X,x}^\times$, where $\cO_{X,x}$ is the stalk (of monoids) at $x$ and $\cO_{X,x}^\times$ denotes its group of invertible elements. We define the \emph{rank of $X$} as $\rk X:=\min_{x\in\tilde X}\{\rk x\}$ and we define
$$ \tilde X^{\rk} \ := \ \coprod_{\rk x=\rk\tilde X} \Spec_\Mo\cO_{X,x}^\times, $$
which is a sub-$\Mo$-scheme of $\tilde X$ whose underlying topological space is discrete.

\begin{df}
A \emph{strong morphism} $\varphi:\cX\to\cY$ is a pair $(\tilde f, f)$, where $\tilde f:\tilde X^\rk\to\tilde Y^\rk$ is a morphism of $\Mo$-schemes and $f:X\to Y$ is a morphism of schemes such that
$$\xymatrix{\tilde X^\rk_\Z\ar[rr]^{\tilde f_\Z}\ar[d]_{e_X} && \tilde Y^\rk_\Z\ar[d]^{e_Y}\\ X\ar[rr]^f && Y} $$
commutes. We denote the category whose objects are $\F_1$-schemes (locally of finite type) and whose morphisms are strong morphisms by $\Sch^\str_\Fun$.
\end{df}

Recall that $\ast_\Mo=\Spec_\Mo\{0,1\}$ is the terminal object in $\Sch_\Mo$ and $\ast_\Z=\Spec\Z$, the terminal object in $\Sch_\Z$. The $\F_1$-scheme $(\ast_\Mo,\ast_\Z,\id_{\ast_\Z})$ is the terminal object in $\Sch_\Fun^\str$, and we denote it by $\ast_\Fun$ or $\Spec_\Fun\F_1$. 

If $\cX=(\tilde X,X,e_X)$ is an $\F_1$-scheme such that $e_X:\tilde X^\rk_\Z\to X$ is an isomorphism, we say that $\cX$ is \emph{of pure rank}, and we denote the full subcategory of those $\F_1$-schemes in $\Sch_\Fun^\str$ by $\Sch_\Fun^\rk$. If $\cX$ is of pure rank and $(\tilde f,f):\cX\to\cY$ is a strong homomorphism between $\F_1$-schemes, then $f\circ e_X=e_Y\circ\tilde f$. Since $e_X$ is an isomorphism, $f=e_Y\circ\tilde f\circ e_X^{-1}$, and we obtain:

\begin{lemma}\label{strong_by_mo}
Let $\cX=(\tilde X,X,e_X)$ and $\cY=(\tilde Y,Y,e_Y)$ be $\F_1$-schemes and $\cX$ be of pure rank. The map $\Hom_\Fun^\str(\cX,\cY)\to\Hom_\Mo(\tilde X,\tilde Y)$ sending $(\tilde f,f)$ to $\tilde f$ is a bijection. \qed
\end{lemma}

Thus we can also consider $\Sch_\Fun^\rk$ as a full subcategory of $\Sch_\Mo$. Its objects are characterized as those $\Mo$-schemes $\tilde X$ for which there is a number $r$ such that $\tilde X$ is a disjoint union of $\Mo$-schemes of the form $\Spec_\Mo(\{0\}\cup H)$, where $H$ is an abelian group of rank $r$. We define
$$ \cY(\F_1) \ := \ \Hom^{\str}_\Fun(\ast_\Fun,\cY) $$
for $\F_1$-schemes $\cY$, and, more generally, $\cY(H):=\Hom^\str_\Fun\bigl(\Spec_\Mo(\{0\}\cup H),\; \cY\bigr)$ for an abelian group $H$ of finite rank and $\cY(\cX):=\Hom^\str_\Fun(\cX,\cY)$ for every $\F_1$-scheme $\cX$ of pure rank.

\begin{lemma}\label{F_1_points}
 Let $\cY=(\tilde Y,Y,e_Y)$ be an $\F_1$-scheme. Then $\cY(\F_1)$ equals the set of points of $\tilde Y^\rk$.
\end{lemma}

\begin{proof}
The $\F_1$-scheme $\ast_\Mo=\Spec_\Mo\{0,1\}$ has one point, namely the unique prime ideal $\{0\}$ of $\{0,1\}$. For every choice of image $y$ of $\{0\}$ in $\tilde Y^\rk$, the stalk $\cO_{Y,y}$ is of the form $\{0\}\cup H$ for an abelian group $H$ and  there is consequently a unique monoid homomorphism $\cO_{Y,y}\to \{0,1\}$ sending $0$ to $0$ and $H$ to $1$.
\end{proof}

If $\cX$ is defined by a torified variety, then $\tilde X^\rk$ corresponds to the tori of lowest dimension in the torification. For split tori, affine spaces, Grassmannians and split reductive groups we described a torification and their tori of lowest dimension in the examples of the previous section. Thus we obtain a solution to Problem A and Problem B, part (i), from the introduction.

\begin{thm} \label{thmA}
In $\Sch_\Fun^\str$, there are objects $\G^r_{m,\Fun}$, $\A^d_\Fun$ (Example \ref{fun_from_toric}), $\Gr(k,n)_\Fun$ (Example \ref{fun_from_schubert}) and there is an object $\cG$ for every split reductive group $G$ (Example \ref{fun_from_split_red}) such that
 \begin{align*}
  \G^r_{m,\Fun}\otimes_\Fun\Z&\simeq \Gm^r&& \text{and} &\G^r_{m,\Fun}(\F_1)&\simeq\ast,  \vspace{4pt}\\
  \A^d_\Fun\otimes_\Fun\Z&\simeq \A^d&& \text{and} &\A^d_\Fun(\F_1)&\simeq\ast,\vspace{4pt}\\
  \Gr(k,n)_\Fun\otimes_\Fun\Z&\simeq \Gr(k,n)&& \text{and} &\Gr(k,n)_\Fun(\F_1)&\simeq M_{k,n}.\vspace{4pt}\\
  \cG\otimes_\Fun\Z&\simeq G&&\text{and}&\cG(\F_1)&\simeq W\quad \text{ (as sets)}.\qed
 \end{align*}
\end{thm}

\begin{rem}\label{models_of_toric}
 For a toric variety $X$, we have defined two different models $\cX_1=(\tilde X_1,X,e_1)$ (Example \ref{ex:fun_from_mo}) and $\cX_2=(\tilde X_2,X,e_2)$ (Example  \ref{fun_from_toric}) of $X$ as $\F_1$-schemes, where $\tilde X_1$ is a connected $\Mo$-scheme and $\tilde X_2$ is a discrete $\Mo$-scheme. However, both $\tilde X_1^\rk$ and $\tilde X_2^\rk$ are discrete and there is an isomorphism $\tilde i:\tilde X_1^\rk\to\tilde X_2^\rk$ of $\Mo$-scheme such that $(\tilde i,\id_X)$ is a strong morphism. This shows that $\cX_1$ and $\cX_2$ are isomorphic in $\Sch_\Fun^\str$.
\end{rem}


\section{Weak morphisms}
\label{weak_morphisms}

\noindent
In this section, we introduce a second notion of morphism between $\F_1$-schemes, which allows us to define all split reductive groups as group objects over $\F_1$. We start with proving some useful facts.

\begin{lemma}\label{injective}
 Let $(\tilde X,X,e_X)$ be a scheme over $\F_1$. As a map between the underlying topological spaces, $e_X:\tilde X_\Z\to X$ is injective.
\end{lemma}

\begin{proof}
 Assume $e_X(x_1)=e_X(x_2)$ for two points $x_1$ and $x_2$ of $\tilde X_\Z$. Then there is a field $k$ and two morphisms $\Spec k\to\tilde X_\Z$ whose images are $\{x_1\}$ and $\{x_2\}$, respectively. Since $e_X$ induces an isomorphism $\tilde X_\Z(k)\simeq X(k)$, the two morphism must have the same image, and thus $x_1=x_2$.
\end{proof}

Let $X^\rk$ denote the image of $e_X: \tilde X^\rk_\Z\to X$. For every point $x\in\tilde X$, let $\{x\}_\Z$ be the corresponding subscheme of $\tilde X_\Z$. We write $e_X(x)$ for the image $e_X(\{x\}_\Z)$ in $X$. By a theorem of Chevalley, the images of constructible sets are constructible. Since $\{x\}_\Z$ is connected, $e_X(x)$ is connected, too, and thus locally closed. This shows that $e_X(x)$ is a subscheme of $X$.

\begin{lemma}\label{disjoint_union}
The image of $\tilde X^\rk$ under $e_X$ is a disjoint union
$$ X^\rk\ =\ \coprod_{x\in\tilde X^\rk}e_X(x). $$
\end{lemma}

\begin{proof}
Since the rank of a point $x\in\tilde X$ equals the dimension of the subscheme $\{x\}_\Z$ of $\tilde X_\Z$ and $e_X$ is injective by Lemma \ref{injective}, the dimension of $e_X(x)$ equals the rank of $X$ for all $x\in \tilde X^\rk$. Since $e_X(x)$ and $e_X(y)$ are disjoint and of equal dimension for two different points $x,y\in\tilde X^\rk$, their union is not connected. Since $X$ is locally of finite type, the image of $e_G$ is a locally finite disjoint union of subschemes of the form $e_G(x)$ with $x\in\tilde X^\rk$. Thus the lemma follows.
\end{proof}

\begin{rem}
 The previous two lemmas show that $e_X$ is an injective map between the underlying topological spaces of $\tilde X_\Z$ and $X$ whose image is locally closed. For to show that $e_X$ is an immersion, we need to show that all morphism between stalks are surjective. It is not clear to me whether this holds true in general. If it holds true, we can identify $\tilde X^\rk_\Z$ and $X^\rk$ via $e_X$. Further it implies that the different definitions of torified schemes (locally of finite type) given in the present text and in \cite{LL09} coincide, cf.\ Remark \ref{different_defs}. 
\end{rem}

Let $\cX=(\tilde X,X,e_X)$ and $\cY=(\tilde Y,Y,e_Y)$ be $\F_1$-schemes. The unique morphism $\Spec_\Mo\cO^\times_{X,x}\to\ast_{\Mo}$ induces a morphism 
$$ \tilde X^\rk\ =\ \coprod_{x\in\tilde X^\rk}\Spec_\Mo\cO^\times_{X,x}\ \longrightarrow \ \ast_\cX \ := \ \coprod_{x\in\tilde X^\rk}\ast_{\Mo}. $$
Given $\tilde f:\tilde X^\rk\to\tilde Y^\rk$, there is a unique morphism $\ast_\cX\to\ast_\cY$ such that
$$\xymatrix@R=1,2pc{\tilde X^\rk\ar[rr]^{\tilde f}\ar[d] && \tilde Y^\rk\ar[d]\\ \ast_\cX\ar[rr] && \ast_\cY} $$
commutes. 

The unique morphism $e_X(x)\to\ast_{\Z}$ to the terminal object $\ast_{\Z}=\Spec\Z$ in $\Sch_\Z$ induces a morphism 
$$ X^\rk\ =\ \coprod_{x\in\tilde X^\rk}e_X(x)\ \longrightarrow \ (\ast_\cX)_\Z \ = \ \coprod_{x\in\tilde X^\rk}\ast_{\Z}. $$

\begin{df}
A \emph{weak morphism} $\varphi:\cX\to\cY$ is a pair $\varphi=(\tilde f, f)$, where $\tilde f:\tilde X^\rk\to\tilde Y^\rk$ is a morphism of $\Mo$-schemes and $f:X\to Y$ is a morphism of schemes that restricts to a morphism $f:X^\rk\to Y^\rk$ such that
$$\xymatrix@C=3pc@R=1pc{\tilde X^\rk_\Z\ar[rr]^{\tilde f_\Z}\ar[dr] && \tilde Y^\rk_\Z\ar[dr]\\
                 &(\ast_\cX)_\Z\ar[rr] && (\ast_\cY)_\Z \\  X^\rk\ar[rr]^f\ar[ur] && Y^\rk\ar[ur]} $$
commutes. We denote by $\Sch_\Fun^\weak$ the category whose objects are $\F_1$-schemes (locally of finite type) and whose morphisms are weak morphisms.
\end{df}


\section{Cartesian categories}
\label{cartesian_categories}

\noindent
We reason that the categories we invented are cartesian and introduce certain cartesian functors that allow us to pass from group objects from one category to group objects of another category by means of Proposition \ref{cart_functor}. The category $\Sch_\Z$ has finite products and $\ast_\Z=\Spec\Z$ as  terminal object. Thus $\Sch_\Z$ is cartesian. We already reasoned in section \ref{f1-schemes} that $\Sch_\Mo$ is cartesian and that the base extension functor $-\otimes_\Mo\Z:\Sch_\Mo\to\Sch_\Z$ is cartesian.

Since the evaluation is an isomorphism for every $\Fun$-scheme of pure rank, the product in $\Sch_\Fun^\rk$ is given by 
$$ (\tilde X,X,e_X) \ \times \ (\tilde Y,Y,e_Y) \ = \ (\tilde X\times\tilde Y,\;X\times Y,\;(e_X,e_Y)). $$  
Since $(\tilde X\times\tilde Y)^\rk=\tilde X^\rk\times\tilde Y^\rk$, the products in $\Sch_\Fun^\str$ and $\Sch_\Fun^\weak$ are realized by the same formula. The terminal object in all three categories is $\ast_\Fun$. It follows that the categories $\Sch_\Fun^\rk$, $\Sch_\Fun^\str$ and $\Sch_\Fun^\weak$ are cartesian and that the inclusion functor $\Sch_\Fun^\rk\hookrightarrow\Sch_\Fun^\str$ is cartesian. Every strong morphism is a weak morphism, thus $\Sch_\Fun^\str$ is a subcategory of $\Sch_\Fun^\weak$. Consequently, the inclusion functor $\Sch_\Fun^\str\hookrightarrow\Sch_\Fun^\weak$ is cartesian.

Recall that we defined the base extension of an $\F_1$-scheme $\cX=(\tilde X,X,e_X)$ as $\cX_\Z=X$. To extend this to a functor $-\otimes_\Fun\Z:\Sch_\Fun^\weak\to\Sch_\Z$, we define the base extension of a weak morphism $\varphi=(\tilde f,f): \cX\to\cY$ to $\Z$ as $\varphi_\Z:=f:\cX_\Z\to\cY_\Z$. This yields a cartesian functor.

We introduce a functor $(-)^\rk: \Sch_\Fun^\weak\to\Sch_\Fun^\rk$ that associates to an $\F_1$-scheme $\cX=(\tilde X,X,e_X)$, the $\F_1$-scheme $\cX^\rk=(\tilde X^\rk,\tilde X^\rk_\Z,\id)$ of pure rank and to a weak morphism $\varphi=(\tilde f, f)$ the strong morphism $\varphi^\rk=(\tilde f,\tilde f_\Z)$. This functor is cartesian.

We subsume these cartesian functors in the following diagram:
$$ \xymatrix@C=3pc{\Sch_\Fun^\rk\ar@{}|\subset[r] &\Sch_\Fun^\str\ar@{}|\subset[r]  &\Sch_\Fun^\weak\ar[rr]^{-\otimes_\Fun\Z}\ar@/_2pc/[ll]_{(-)^\rk}  &&\Sch_\Z. } $$

Note that the composition of cartesian functors is cartesian. Thus the base extension to $\Z$ restricted to $\Sch_\Fun^\rk$ or $\Sch_\Fun^\str$ is cartesian, too. A consequence of the fact that every strong morphism $\varphi:\cX\to\cY$ factors through $\cY^\rk$ when $\cX$ is of pure rank is that $\Hom_\Fun^\rk(\cX,\cY^\rk)=\Hom_\Fun^\str(\cX,\cY)$.

An immediate consequence of Proposition \ref{cart_functor} is the following key property that will be of important (implicit) use for the theory of algebraic groups over $\F_1$ as introduced in the next section.

\begin{lemma}
Let $\cG$ be an algebraic group over $\F_1$ with group law $\mu$. Then $\cG_\Z$ is an algebraic group (over $\Z$) with group law $\mu_\Z$, and $\cG^\rk$ is a group object in $\Sch_\Fun^\rk$ with group law $\mu^\rk$. In particular, for every $\F_1$-scheme $\cX$ of pure rank, $\cG(\cX)=\Hom_\Fun^\rk(\cX,\cG^\rk)$ inherits the structure of a group.\qed
\end{lemma}

\section{Algebraic groups over $\F_1$}
\label{algebraic_f1_groups}

\noindent
The subject of this section is to establish the notion and various examples of algebraic groups over $\F_1$. 

\begin{df}
 A \emph{group scheme over $\Fun$} is a group object in $\Sch_\Fun^\weak$. Let $G$ be a group scheme with group law $m$. If there is a group scheme $(\cG,\mu)$ over $\Fun$ such that $\cG_\Z\simeq G$ as group schemes (over $\Z$), then we say that $\cG$ is a \emph{model of $G$ over $\F_1$}. If $\mu$ is a strong morphism, we say that $\cG$ is a \emph{canonical model of $G$ over $\F_1$}. A group scheme $\cG$ over $\F_1$ is called an \emph{algebraic group over $\F_1$} if $\cG_\Z$ is an algebraic group (over $\Z$).
\end{df}

For a group $W$, denote by $W_\Z$ the constant group scheme of $W$, i.e.\ the scheme $W_\Z=\coprod_W\Spec\Z$ together with the obvious group law. Clearly, we have a model of the constant group scheme in $\Sch_\Fun^\str$. More precisely:

\begin{lemma}\label{const_group}
For every group $W$, there is a group law $\mu: W_\Fun\times W_\Fun\to W_\Fun$ in $\Sch_\Fun^\str$ for the $\F_1$-scheme $W_\Fun=\coprod_W\ast_\Fun$ such that $W_\Fun\otimes_\Fun\Z\simeq W_\Z$ as group schemes and $W_\Fun(\Fun)\simeq W$ as groups.\qed
\end{lemma}

\begin{lemma}\label{split_tori}
For every $r\geq0$, there is a group law $\G_{m,\Fun}^r\times\G_{m,\Fun}^r\to\G_{m,\Fun}^r$ in $\Sch_\Fun^\str$ such that $\G_{m,\Fun}^r\otimes_\Fun\Z\simeq\Gm^r$ as algebraic groups. The group $\G_{m,\Fun}^r(\Fun)$ is the trivial group.
\end{lemma}

\begin{proof}
 Since $\G_{m,\Mo}^r=\Spec_\Mo(\{0\}\cup H)$, where $H=\{T_1^{n_1}\cdot\dotsb\cdot T_r^{n_r}\mid n_1,\dotsc,n_r\in\Z\}$ is the free abelian group in $r$ generators, it has precisely one point, namely the unique prime ideal $\{0\}$ of $H$. The stalk $\cO_{\G_{m,\Mo}^r,\{0\}}$ is equal to $H$. If we let $m$ be the group law of $\Gm^r$ (in $\Sch_\Z$) and define $\tilde m$ topologically as the trivial map and on the stalk as the map $\tilde m ^\#: H\to H\wedge H$ sending an element $h\in H$ to $(h,h)\in H$ and $0$ to $0$, then
$$ \xymatrix@R=1,5pc{(\G_{m,\Mo}^r)_\Z\times(\G_{m,\Mo}^r)_\Z\ar[rr]^{\tilde m_\Z}\ar@{=}[d] && (\G_{m,\Mo}^r)_\Z\ar@{=}[d] \\ \Gm^r\times\Gm^r\ar[rr]^m && \Gm^r} $$
commutes. If $\epsilon:\ast_\Z\to\Gm^r$ and $\iota:\Gm^r\to\Gm^r$ are the morphisms of Proposition \ref{alt_group}, then it is easily seen that they extend to morphisms $(\tilde\epsilon,\epsilon):\ast_\Fun\to\G_{m,\Fun}^r$ and $(\tilde\iota,\iota):\G_{m,\Fun}^r\to\G_{m,\Fun}^r$ of $\F_1$-schemes that satisfy the definition of an algebraic group. Thus $\G_{m,\Fun}^r$ is a group object in $\Sch_\Fun^\str$ and $\G_{m,\Fun}^r\otimes_\Fun\Z=\Gm^r$ as algebraic groups by construction.

Since $\G_{m,\Mo}^r$ has only one point and since there is only one homomorphism of monoids $\{0\}\amalg H\to\{0,1\}$, $\G_{m,\Fun}^r(\Fun)$ is the trivial group.
\end{proof}

This proof generalizes to show the more general lemma.

\begin{lemma}
 Let $H$ be a finitely generated abelian group. Then $\cG:=(\Spec_\Mo(\{0\}\amalg H), \Spec \Z[H],\id)$ is a group scheme over $\F_1$ with $\cG_\Z=\Spec\Z[H]$ and $\cG(\F_1)$ being the trivial group.\qed
\end{lemma}

Let $G$ be a group scheme with group law $m$. If $e_G: \coprod_{i\in I}\Gm^{d_i}\to G$ is a torification of $G$, then let $(\tilde G,G,e_G)$ be the associated $\F_1$-scheme (cf.\ Lemma \ref{fun_from_torified}). Let $N=G^\rk$ be the image of $\tilde G^\rk_\Z$ under $e_G$, put $r=\min_{i\in I}\{d_i\}$ and $I^\rk=\{i\in I\mid d_i=r\}$, and denote by $e_G^\rk:\tilde G^\rk_\Z\to N$ the restriction of $e_G$ to $\tilde G^\rk_\Z=\coprod_{i\in I^\rk}\Gm^r$. If $e_G^\rk$ is an isomorphism of schemes and $m$ restricts to a group law of $N$, then define $W:=N/T$, where $T=N^0\simeq\Gm^r$ is the connected component of $N\simeq \coprod_{i\in I^\rk}\Gm^r$. Thus we can identify $I^\rk$ with $W$ and write $N\simeq\coprod_W\Gm^r$. 

\begin{thm}\label{torified_to_alg_groups}
Let $(G,m)$ be a group scheme with torification $e_G$ such that $m$ restricts to a group law of $N=G^\rk$. Put $T=N^0$, $W=N/T$ and $r=\dim T$. Assume that $e_G^\rk:\coprod_W\Gm^r\to N$ is an isomorphism.
\begin{enumerate} 
 \item\label{tor1} There is a model $\cG$ of $G$ over $\F_1$ and $\cG(\F_1)\simeq W$.
 \item\label{tor2} If $N\simeq T\rtimes_\theta W$ for some $\theta: W\times T\to T$, then there is a canonical model $\cG$ of $G$ over $\F_1$. If $e_G'$ is another torification of $G$ that satisfies the hypotheses for to have a canonical model $\cG'$ over $\F_1$ and if there is an group automorphism $f:G\to G$ such that $f(N)=N'$, where $N'=G^\rk$ w.r.t.\ $e_G'$, then $\cG$ and $\cG'$ are isomorphic in $\Sch_\Fun^str$.
\end{enumerate} 
\end{thm}

\begin{proof}
We begin with \eqref{tor1}. Let $\cG=(\tilde G,G,e_G)$ be the $\F_1$-scheme associated to the torification $e_G$. Lemma \ref{const_group} provides a canonical model $W_\Fun=(W_\Mo,W_\Z,e_W)$ of $W_\Z$ over $\F_1$. In particular, we have a group object $W_\Mo=\coprod_W\ast_\Mo$ with group law $\tilde m_W$ in $\Sch_\Mo$. Lemma \ref{split_tori} provides a canonical model $\cT=(\tilde T,T,e_T)$ of $T\simeq\Gm^r$ over $\F_1$ and in particular a group object $\tilde T\simeq\G_{m,\Mo}^r$ with $\tilde m_T$ in $\Sch_\Mo$.

Choose any group action $\tilde\theta: W_\Mo\times\tilde T\to\tilde T$ that respects $\tilde m_T$, e.g.\ the projection to the second component, which is always possible. Then $\tilde\theta$ defines the semidirect product $\tilde N:=\tilde T\rtimes_{\tilde\theta}W_\Mo$ with group law $\tilde m$ as a group object in $\Sch_\Mo$. The diagram
$$\xymatrix@R=1pc{\tilde N_\Z\times\tilde N_\Z\ar[rr]^{\tilde m_\Z}\ar[dr] && \tilde N_\Z\ar[dr]\\
                 &W_\Z\times W_\Z\ar[rr] && W_\Z \\  N\times N\ar[rr]^m\ar[ur] && N^\rk\ar[ur]} $$
commutes. This shows that $\mu=(\tilde m,m): \cG\times \cG\to\cG$ is a weak morphism.

To verify that $\mu$ is a group law, let $\epsilon$ and $\iota$ be the unit and the inversion of $G$ and let $\tilde\epsilon$ and $\tilde\iota$ be the unit and the inversion of $\tilde G$ (cf.\ Proposition \ref{alt_group}). It is easily seen that $(\tilde\epsilon,\epsilon)$ and $(\tilde\iota,\iota)$ are weak morphisms (the former one being even a strong morphism) and verify the conditions of Proposition \ref{alt_group}. Thus $\cG$ is a group scheme over $\F_1$ whose base extension to $\Z$ is $G$.

By Lemma \ref{F_1_points}, $\cN(\F_1)=\tilde N^\rk=W$ as sets. Since every strong morphism $\ast_\Fun\to\G_{m,\Fun}^r$ is trivial, a strong morphism $\ast_\Fun\to\cN$ factorizes through $\ast_\Mo\to W_\Fun$, thus $\cN(\F_1)=W$ as groups.

We proceed with \eqref{tor2}. We need to define $\tilde\theta$ as above such that $(\tilde\theta,\theta)$ is a strong morphism.

Fix an isomorphism $T\simeq\Gm^r$. For every $w\in W$, the restriction $\theta_w:\Gm^r\to\Gm^r$ of $\theta$ to the component of $\coprod_W\Gm^r$ corresponding to $w$ is a homomorphism of group schemes, since $\theta$ respects the group law of $\Gm^r$. Let $H=\{T_1^{n_1}\cdot\dotsb\cdot T_r^{n_r}\mid n_1,\dotsc,n_r\in\Z\}$, then $\theta_w$ induces an automorphism $\theta_w^\#$ of the Hopf algebra $\Z[H]$, which restricts to a group automorphism $\tilde\theta_w^\#:H\to H$ of the group-like elements $H$ of $\Z[H]$. This defines a homomorphism $\tilde\theta_w: \G_{m,\Mo}^r\to\G_{m,\Mo}^r$ of group objects in $\Sch_\Mo$ that base extends to $(\tilde\theta_w)_\Z=\theta_w$. We obtain a morphism of $\Mo$-schemes 
$$ \tilde\theta \ = \ \coprod_{w\in W}\tilde\theta_w: \ W_\Mo\times\G_{m,\Mo}^r \ = \ \coprod_W \G_{m,\Mo}^r \to \G_{m,\Mo}^r $$
whose base extension to $\Z$ is $\tilde\theta_\Z=\theta$. Therefore, $(\tilde\theta,\theta):W_\Fun\times\G_{m,\Fun}^r\to\G_{m,\Fun}^r$ is a strong morphism. Note that $\theta$ is a group action, and it respects the group law of $\G_{m,\Mo}^r$, since the restrictions $\tilde\theta_w$ are homomorphisms of group objects for all $w\in W$. 

That canonical models associated to $e_G$ and $e_G'$ as in the theorem are isomorphic is reasoned as follows. Since $(\tilde\theta,\theta)$ is a strong morphism, we have that $\theta\circ e_N=e_T\circ\tilde\theta$, where $e_N=(e_W,e_T)$. But $e_T$ is an isomorphism, thus $\tilde\theta_\Z$ is determined by $\theta$. Since the automorphism $f:G\to G$ restricts to an isomorphism between $N$ and $N'$, which is an isomorphism of algebraic groups. Going through the construction of $\tilde\theta$, we see that this defines already a morphism $\tilde f:\tilde N\to \tilde N'$ of $\Mo$-schemes such that $(\tilde f,f):\cG\to\cG'$ is a strong morphism that is an isomorphism of group schemes over $\F_1$.
\end{proof}

We have some immediate corollaries.

\begin{cor}
 If $N$ is an extension of a constant group scheme $W_\Z$ associated to a group $W$ by a split torus, then $N$ has a model $\cN$ over $\F_1$. The group $\cN(\F_1)$ is isomorphic to $W$. If $N$ is a split extension, it has a canonical model over $\F_1$.\qed
\end{cor}

\begin{rem}\label{semidirect_remark}
Let $N=\Gm^r\rtimes_\theta W_\Z$ be as in the theorem. Let $\cN=(\tilde N,N,e_N)$ be the canonical model of $N$ over $\F_1$. Since $\tilde N^\rk=N$, the group law $\mu$ extends (trivially) to a group law of $\cN$ in the category $\Sch_\Fun^\nat$ whose objects are $\F_1$-schemes and whose morphisms are pairs $(\tilde f,f):(\tilde X,X,e_X)\to(\tilde Y,Y,e_Y)$, where $\tilde f:\tilde X\to\tilde Y$ is a morphism between $\Mo$-schemes and $f:X\to Y$ is a morphism between schemes such that $e_Y\circ\tilde f_\Z=f\circ e_X$. More generally, we can substitute $\Gm^r$ by a group scheme of the form $\Spec\Z[H]$, where $H$ is a finitely generated abelian group.

However, semidirect products of a group schemes of the form $\Spec\Z[H]$ with a finite constant group scheme seem to be the only algebraic groups that allow models in $\Sch_\Fun^\nat$, cf.\ the explanations in \cite[section 6.1]{LL09}. The following implications of Theorem \ref{torified_to_alg_groups} show that $\Sch_\Fun^\str$ allows a larger class of group objects.
\end{rem}

Let $\Ga$ be the additive group scheme. We say that $G$ is a \emph{successive extension of additive groups} if there is a sequence of subgroups $0=G_0<G_1<\dotsb<G_n=G $ such that $G_{i-1}$ is a normal subgroup of $G_i$ and $G_i/G_{i-1}\simeq\Ga$ for all $i=1,\dotsc,n$. 

\begin{cor}\label{additive}
 Let $G$ be an algebraic group that is a successive extension of additive groups. Then it has a canonical model $\cG$ over $\F_1$. The group $\cG(\F_1)$ is the trivial group.
\end{cor}

\begin{proof}
As a variety, $G$ is isomorphic to $\A^n$, where $n$ is the number of subgroups in the filtration of $G$. Let $\epsilon:\ast_\Z\to G$ be the unit of $G$. We can choose a torification of $G\simeq\A^n$ such that $\epsilon$ is the unique $0$-dimensional torus of the torification. This defines an $\F_1$-scheme $\cG=(\tilde G,G,e_G)$ such that $\tilde G^\rk\simeq\ast_\Mo$, and $\tilde G^\rk_\Z\simeq\ast_\Z\stackrel\epsilon\to G$ equals $e_G^\rk$. Thus Theorem \ref{torified_to_alg_groups} \eqref{tor2} applies and proves that $G$ has a canonical model over $\F_1$. Clearly, $\cG(\Fun)=\Hom(\ast_\Mo,\tilde G^\rk)$ is the trivial group.
\end{proof}

If $G$ is a reductive group with Weyl group $W=N(\Z)/T(\Z)$, then we say that \emph{$W$ lifts along $\sigma$} if the exact sequence of groups
\begin{equation*}\xymatrix@C=4pc{1\ar[r]&T\ar[r]&N\ar[r]&W_\Z\ar[r]\ar@{-->}@/_1pc/[l]_\sigma&1.} \end{equation*}
splits. As a consequence of this theorem we obtain the following solution to Problem B of the introduction.

\begin{thm}\label{thmB}
Let $G$ be a split reductive group with group law $m:G\times G\to G$ and Weyl group $W$. Then $G$ has a model $\cG=(\tilde G,G,e_G)$ over $\F_1$ and there is an isomorphism $\cG(\F_1)\simeq W$ of groups such that 
$$ \begin{array}{cccccccc}\hspace{0pt}\sigma: & N(\Z)/T(\Z) \ = \ W & \stackrel\sim\longrightarrow & \cG(\F_1) & \stackrel{-\otimes_\Fun\Z}\longrightarrow& \cG_\Z(\Z)&=& G(\Z) \\  &&&  (\ast_\Fun\to\cG)& \longmapsto & (\ast_\Z\to\cG_\Z) \end{array} $$
maps each coset $nT(\Z)$ in $N(\Z)/T(\Z)$ to an element of $nT(\Z)\subset G(\Z)$.

If the Weyl group lifts along a group homomorphism $\sigma':W\to N(\Z)$, then there is a canonical model $\cG$ of $G$ over $\F_1$ and an isomorphism $\cG(\F_1)\simeq W$ such that $\sigma$ coincides with $\sigma'$.
\end{thm}

\begin{proof}
 In Example \ref{fun_from_split_red}, we endowed a split reductive group $G$ with a torification $e_G$ that restricts to a torification $e_N$ of the normalizer $N$ of a maximal split torus $T$. Since $G^\rk=N$ w.r.t.\ $e_G$, we have that $e_N=e_G^\rk$, and  $e_G^\rk: \coprod_W\Gm^rt\to N$ is an isomorphism. Theorem \ref{torified_to_alg_groups} \eqref{tor1} shows that $G$ has a model $\cG$ over $\F_1$ and that $\cG(\F_1)=\tilde G^\rk=W$. That $\sigma$ maps each coset $nT(\Z)$ in $N(\Z)/T(\Z)$ to an element of $nT(\Z)\subset G(\Z)$ is clear by the construction of $(\cG,\mu)$ and $\cG(\F_1)\simeq W$ in the proof of Theorem \ref{torified_to_alg_groups}.
 
If the Weyl group lifts along a group homomorphism $\sigma':W\to N(\Z)$, then $W$ can be considered as a subgroup of $N(\Z)$, or, equivalently, $W_\Z$ can be considered as a subgroup of $N$. Since $T$ is normal in $N$, $W_\Z$ acts by conjugation on $T$. The conjugation $\theta$ respects the group law of $T$. Thus $N=T\rtimes_\theta W_\Z$ and we can apply Theorem \ref{torified_to_alg_groups} \eqref{tor2} to obtain a canonical model. Again, it is clear from the proof of Theorem \ref{torified_to_alg_groups} that $\sigma$ and $\sigma'$ coincide.
\end{proof}

\section{Parabolic subgroups of $\GL(n)$}
\label{parabolic_subgroups}

\noindent
In this last section, we will investigate Problem C from the introduction. We show that it can be solved within the framework of this paper.

\begin{lemma}\label{parabolic}
Let $P$ be a parabolic subgroup of $\GL(n)$ of type $(k_1,\dotsc,k_r)$. Then $P$ has a canonical model $\cP(\F_1)$ over $\F_1$ and $\cP(\F_1)\simeq S_{k_1}\times\dotsb\times S_{k_r}$. In particular, $\GL(n)$ has a canonical model $\cG$ over $\F_1$ and $\cG(\Fun)\simeq S_n$.
\end{lemma}

\begin{proof}
A parabolic subgroup $P$ of $\GL(n)$ of type $(k_1,\dotsc,k_r)$ is an extension of $M=\GL(k_1)\times \dotsb\times\GL(k_r)$ by a successive extension $U$ of additive groups. Hence, the parabolic subgroup $P$ has a maximal split torus $T$ of rank $n=k_1+\dotsb+k_r$. Let $N=\coprod_W T$ be the normalizer of $T$. Then $W=N(\Z)/T(\Z)\simeq S_{k_1}\times\dotsb\times S_{k_r}$ and the sequence $1\to T(\Z)\to N(\Z)\to W\to1$ splits.

As a variety, $P\simeq M\times U$, thus the product torification of torifications of $M$ and $U$ is a torification of $P$. Choose a torification of $M$ relative to the torus $T$ as described in Example \ref{fun_from_split_red} and for $U$ as described in the proof of Corollary \ref{additive}. Then the product torification $e_P$ of $P$ defines an $F_1$-scheme $\cP=(\tilde P,P,e_P)$ such that $e_P$ restricts to an isomorphism $\tilde P^\rk_\Z\simeq N$. Thus Theorem \ref{torified_to_alg_groups} \eqref{tor2} applies and implies the statement of the proposition. (Note that $\GL(n)$ is a parabolic subgroup of type $(n)$ of $\GL(n)$).
\end{proof}

Recall from Example \ref{fun_from_split_red} that a choice of a maximal split torus $T$ in $\GL(n)$ and a Borel subgroup $B$ containing $T$ leads to a Bruhat decomposition $\coprod_W BwB\to\GL(n)$, where $W\simeq S_n$ is the Weyl group of $\GL(n)$. This leads further to a torification $e_G$ of $\GL(n)$ and defines an $\F_1$-scheme $\cG=(\tilde G,G,e_G)$. By Theorem \ref{thmB}, there is a group law $\mu=(\tilde m,m)$ of $\cG$ such that $\cG$ is a canonical model of $G$. This canonical model depends a priori on the choice of $T$ and $B$, but since all maximal split tori in $\GL(n)$ are conjugated, Theorem \ref{torified_to_alg_groups} \eqref{tor2} implies that the canonical model $\cG$ is unique up to isomorphism. Let $P$ be a parabolic subgroup of type $(k,n-k)$ of $\GL(n)$ that contains $T$ and $B$, and let $\cP=(\tilde P,P,e_P)$ be the canonical model as described in Lemma \ref{parabolic}. Let $l: P\times G\to G$ be the restriction of $m:G\times G\to G$ to the natural action of $P$ on $G$ by left multiplication.


\begin{thm}\label{thmC}
 In the situation as above, the following holds true. 
 \begin{enumerate}
 \item \label{pg1} Then there is a morphism $\tilde l:\tilde P^\rk\times\tilde G^\rk\to\tilde G^\rk$ of $\Mo$-schemes such that $\lambda=(\tilde l,l):\cP\times \cG\to\cG$ is a group action in $\Sch_\Fun^\str$. Taking $\F_1$-points is compatible with the natural group action 
 $$\lambda(\F_1): (S_k\times S_{n-k})\times S_n\to S_n.$$
 \item \label{pg2} There is an $\F_1$-scheme $\cQ$ that is a quotient of $\lambda$. Consequently, 
 $$ \cQ_\Z\simeq\Gr(k,n) \hspace{1cm}\text{and}\hspace{1cm} \cQ(\F_1)\simeq M_{k,n}. $$
 The natural action $\tau:\cG\times\cQ\to\cQ$ on the quotient is compatible with the natural action $\GL(n)\times\Gr(k,n)\to\Gr(k,n)$ and taking $\F_1$-points of $\tau$ is compatible with the natural action 
$$ \tau(\F_1): S_n\times M_{k,n}\to M_{k,n} $$ 
induced by permuting the elements of $M_n=\{1,\ldots,n\}$.
 \end{enumerate}
\end{thm}

\begin{proof}
We begin with \eqref{pg1}. The maximal split torus $T$ is a subgroup of both $P$ and $G$. Its normalizer $N_P$ in $P$ is a subgroup of its normalizer $N$ in $G$. By construction of $\cP=(\tilde P,P,e_P)$, we have that $N_P=\tilde P^\rk_\Z$ (cf.\ Lemma \ref{parabolic}) and by construction of $\cG=(\tilde G,G,e_G)$, we have that $N=\tilde G^\rk_\Z$ (cf.\ Theorem \ref{thmB}). Put $W_{P,\Z}=N_P/T$ and $W_\Z=N/T$. Then we obtain an inclusion $W_{P,\Z}\hookrightarrow W_\Z$ of groups. Since $W_\Z$ lifts to a subgroup of $G$, $N$ is a semidirect product $T\rtimes_{\theta}W_\Z$ along a group action $\theta:W_\Z\times T\to T$. If $\theta_P:W_{P,\Z}\times T\to T$ is the restriction of $\theta$, then $N_P$ is the semidirect product $T\rtimes_{\theta_P}W_{P,\Z}$.
 
These semidirect products define group laws $\tilde m_P$ and $\tilde m$ on $\tilde P^\rk$ and $\tilde G^\rk$, respectively, such that $\tilde P^\rk$ is a subgroup of $\tilde G^\rk$. Consequently, the restriction of $\tilde m$ defines an action $\tilde l:\tilde P^\rk\times\tilde G^\rk\to\tilde G^\rk$. Since $(\tilde m,m)$ is a strong morphism, $\lambda=(\tilde l,l)$ is a strong morphism, too. By Theorem \ref{thmB} and Lemma \ref{parabolic}, taking $\F_1$-points yields $\lambda(\F_1): (S_k\times S_{n-k})\times S_n\to S_n$ as desired.

We proceed with \eqref{pg2}. We construct $\cQ=(\tilde Q,Q,e_Q)$ as follows. Define $Q=\Gr(k,n)$. We review the Schubert decomposition in detail. We have the decompositions 
$$  \coprod_{w\in W_P}BwB\ \longrightarrow \ P\hspace{1,5cm} \text{and} \hspace{1,5cm} \coprod_{w\in W}BwB\ \longrightarrow \ G, $$
where $W_P=W_{P,\Z}(\Z)$, $W=W_\Z(\Z)$ and $w\in W$ is identified with the image of the corresponding point of $W_\Z$ in $G$. These decompositions yield a decomposition
$$ \coprod_{w\in W/W_P} (BwB)\;/\;(BW_{P,\Z}B) \ \longrightarrow \ \Gr(k,n) \ = \ G/P.   $$
The quotients $(BwB)/(BW_{P,\Z}B)$ are affine spaces $\A^{d_w}$ of a certain dimension $d_w$ for every coset $w\in W/W_P$. We obtain a Schubert decomposition of $\Gr(k,n)$ and we refine this decomposition to a torification $e_Q$ whose $0$-dimensional tori coincide with the morphisms $\Gm^0=T/T\hookrightarrow(BwB)/(BW_{P,\Z}B)$ for every $w\in W/W_P$. This torification defines an $\F_1$-scheme $\cQ=(\tilde Q,\Gr(k,n),e_Q)$.

Since the tori of lowest dimension in the torification of $G$ are the immersions $T\hookrightarrow BwB$ for every $w\in W$ and the tori of lowest dimension in the torification of $P$ are the immersions $T\hookrightarrow BwB$ for every $w\in W_P$, the $\Mo$-scheme $\tilde Q^\rk$ is the quotient of the action $\tilde l: \tilde P^\rk\times\tilde G^\rk\to\tilde G^\rk$. Thus $\cQ$ is a quotient of $\lambda$. 

By construction, we have $\cQ_\Z\simeq\Gr(k,n)$. The group $W$ is the Weyl group of $\GL(n)$ and thus naturally isomorphic to $S_n$, and $W_P$ is naturally isomorphic to $S_k\times S_{n-k}$ by Lemma \ref{parabolic}. Thus we have 
$$ \cQ(\F_1) \ \simeq \ \cG(\F_1)/\cP(\F_1) \ \simeq \ W/W_P \ \simeq \ S_n/(S_k\times S_{n-k}). $$

By construction, the natural action $\cG\times\cQ\to\cQ$ is after base extension to $\Z$ compatible with the natural action $G\times Q\to Q$. The identification $M_{k,n}=S_n/(S_k\times S_{n-k})$ yields the natural action of $\cG(\Fun)=S_n$ on $\cQ(\Fun)=M_{k,n}$. 
\end{proof}

\bibliographystyle{plain}

\end{document}